\documentclass[a4paper,11pt]{article}
\usepackage[OT1,T1]{fontenc}
\usepackage{a4wide}
\usepackage{amsmath,amssymb}
\usepackage{mathrsfs}
\usepackage{newcent,fouriernc}
\usepackage{tikz}
\usepackage[pdfborder={0 0 0},pdfpagemode=None,pdfstartview=FitH]{hyperref}

\makeatletter
\renewcommand\section{\@startsection{section}{1}{\z@}
                                  {-3.5ex\@plus-1ex\@minus-.2ex}
                                  {2.3ex\@plus.2ex}
                                  {\normalfont\large\scshape}}
\renewcommand\subsection{\@startsection{subsection}{2}{\z@}
                                     {3.25ex\@plus1ex\@minus.2ex}
                                     {-1em}
                                     {\normalfont\normalsize\bfseries}}
\makeatother

\renewenvironment{description}{\list{}
  {\setlength\leftmargin{2.5em}\setlength\labelwidth{0pt}
  }}
  {\endlist}

\DeclareMathOperator\Hom{Hom}
\DeclareMathOperator\coker{coker}
\DeclareMathOperator\im{im}
\DeclareMathOperator\wt{wt}
\DeclareMathOperator\ad{ad}
\DeclareMathOperator\Irr{Irr}
\DeclareMathOperator\linspan{span}
\newcommand\GL{{\mathbf{GL}}}
\newcommand\dimvec{{\,\underline\dim\,}}
\newcommand\out{\mathrm{out}}
\newcommand\iin{\mathrm{in}}
\newcommand\mmod{\text{\upshape-}\mathrm{mod}}

\begin{document}
\title{Weyl group action and semicanonical bases}
\author{Pierre Baumann}
\date{}
\maketitle

\begin{abstract}
\noindent
Let $\mathbf U$ be the enveloping algebra of a symmetric Kac-Moody
algebra. The Weyl group acts on $\mathbf U$, up to a sign. In addition,
the positive subalgebra $\mathbf U^+$ contains a so-called semicanonical
basis, with remarkable properties. The aim of this paper is to show
that these two structures are as compatible as possible.
\end{abstract}

\section{Introduction}
\label{se:Intro}
\subsection{}
\label{ss:MainResult}
Let $A=(a_{ij})$ be a generalized Cartan matrix, with lines and columns
indexed by a set~$I$. From the datum of $A$, one builds a Kac-Moody
algebra $\mathfrak g$; its derived algebra is generated by
$\mathfrak{sl}_2$-triples $(e_i,h_i,f_i)$, for $i\in I$. Let $\mathbf U$
be the enveloping algebra of $\mathfrak g$ and let $\mathbf U^+$ be the
subalgebra of $\mathbf U$ generated by the elements $e_i$.

Let us fix $i\in I$. The element $s_i=\exp(-e_i)\exp(f_i)\exp(-e_i)$
belongs to the Kac-Moody group attached to $\mathfrak g$. Through the
adjoint action, $s_i$ defines an automorphism $T_i$ of $\mathbf U$.
We set $\mathbf U^+_i=\mathbf U^+\cap T_i^{-1}(\mathbf U^+)$ and
${}_i\mathbf U^+=\mathbf U^+\cap T_i(\mathbf U^+)$; thus $T_i$
restricts to an isomorphism $\mathbf U^+_i\to\,{}_i\mathbf U^+$.
Moreover, the decompositions
$$\mathbf U^+=\mathbf U^+e_i\oplus\mathbf U^+_i
=e_i\mathbf U^+\oplus\;{}_i\mathbf U^+$$
yield projections $\pi_i:\mathbf U^+\to\mathbf U^+_i$ and
${}_i\pi:\mathbf U^+\to\,{}_i\mathbf U^+$.

The algebra $\mathbf U^+$ contains a special basis $\mathbf B$, defined
by Lusztig \cite{Lusztig93} and called the canonical basis. It induces
a basis $\pi_i(\mathbf B)\setminus\{0\}$ in $\mathbf U^+_i$ and a basis
${}_i\pi(\mathbf B)\setminus\{0\}$ in ${}_i\mathbf U^+$. (The two
subspaces $\mathbf U^+_i$ and ${}_i\mathbf U^+$ of $\mathbf U^+$ are in
fact equal, but we distinguish them because the induced bases are
different.) In \cite{Lusztig96}, Lusztig shows that these two bases
correspond to each other under $T_i:\mathbf U^+_i\to\,{}_i\mathbf U^+$.

When $A$ is symmetric, $\mathbf U^+$ can also be endowed with Lusztig's
semicanonical basis \cite{Lusztig00}. Though not being algorithmically
computable, this basis recently attracted some interest because of its
relation with the theory of cluster algebra, see~\cite{Leclerc10} for a
recent survey. The main result of the present paper is a proof that the
above statement about the canonical basis also holds true for the
semicanonical basis.

\subsection{}
\label{ss:Corollaries}
Let $*$ be the antiautomorphism of $\mathbf U^+$ that fixes all the
generators $e_i$. This involution exchanges $(\mathbf U^+e_i,
\mathbf U^+_i,\pi_i)$ and $(e_i\mathbf U^+,{}_i\mathbf U^+,{}_i\pi)$.
It leaves stable the canonical and the semicanonical bases.

Let $B=B(-\infty)$ be the crystal (in the sense of Kashiwara)
associated to the crystal basis of $U^+_q(\mathfrak g)$. Besides
the maps $\wt$, $\varepsilon_i$, $\varphi_i$ and the operations
$\tilde e_i$ and $\tilde f_i$, the set $B$ is endowed with an
involution $b\mapsto b^*$, which reflects the existence of $*$.
This invites us to define the starred operators
$\tilde e_i^*=\bigl(b\mapsto(\tilde e_ib^*)^*\bigr)$ and
$\tilde f_i^*=\bigl(b\mapsto(\tilde f_ib^*)^*\bigr)$. Given
$i\in I$, we set
$$B_i=\{b\in B\mid\varphi_i(b^*)=0\}.$$
In Corollary~3.4.8 in \cite{Saito94}, Saito defines a bijection
$\sigma_i:B_i\to(B_i)^*$ by the rule
$$\sigma_i(b)=\tilde e_i^*{}^{\varepsilon_i(b)}(\tilde f_i)^{\max{}}b.$$

The canonical and semicanonical bases are naturally indexed
by the crystal $B$: to each $b\in B$ correspond elements $G(b)$
and $S(b)$ in the canonical and in the semicanonical bases,
respectively. By Theorem~14.3.2 in~\cite{Lusztig93} and Theorem~3.1
in~\cite{Lusztig00}, both $\{G(b)\mid b\in B\setminus B_i\}$ and
$\{S(b)\mid b\in B\setminus B_i\}$ are bases of $\mathbf U^+e_i$,
therefore both $\{\pi_i(G(b))\mid b\in B_i\}$ and
$\{\pi_i(S(b))\mid b\in B_i\}$ are bases of $\mathbf U^+_i$.
Lusztig's result recalled above can then be given the more precise form:
$$\forall b\in B_i,\qquad T_i(\pi_i(G(b)))=\,{}_i\pi(G(\sigma_ib)).$$
Our aim is to prove the analog formula for the semicanonical basis:
\begin{equation}
\label{eq:LuszThm}
\forall b\in B_i,\qquad T_i(\pi_i(S(b)))=\,{}_i\pi(S(\sigma_ib)).
\end{equation}

We consider the dual framework. Let us denote by $\{S^*(b)\mid b\in B\}$
the dual semicanonical basis, namely the basis of $(\mathbf U^+)^*$ dual
to the semicanonical basis. Then (\ref{eq:LuszThm}) is equivalent to
\begin{equation}
\label{eq:DualForm}
\forall(b',u)\in B_i\times\mathbf U^+_i,\qquad
\langle S(b')^*,u\rangle=\langle S(\sigma_ib')^*,T_i(u)\rangle.
\end{equation}
Taking $u=\pi_i(G(b))$ in (\ref{eq:DualForm}), we get
\begin{equation}
\label{eq:CompBases}
\forall(b,b')\in(B_i)^2,\qquad
\langle S(b')^*,G(b)\rangle=\langle S(\sigma_ib')^*,G(\sigma_ib)\rangle.
\end{equation}

Relation (\ref{eq:CompBases}) constrains the transition matrix between
the canonical and the semicanonical bases.

\section{Adapted filtrations}
\label{se:AdapFilt}
This section is devoted to the statement and to the proof of a key
combinatorial result.

\subsection{}
\label{ss:ConstFun}
Given an algebraic variety $X$ over $\mathbb C$, we denote by $M(X)$
the $\mathbb Q$-vector space consisting of all constructible functions
$f:X\to\mathbb Q$, that is, the $\mathbb Q$-vector space spanned by
the indicator functions of the locally closed subsets of $X$. Let
$\int_X:M(X)\to\mathbb Q$ be the linear form given by
$f\mapsto\sum_{a\in\mathbb Q}a\;\chi(f^{-1}(a))$, where $\chi$ denotes
the Euler characteristic with compact support.

\subsection{}
\label{ss:CombiSetup}
We consider the following setup:
\begin{description}
\item[(a)]
$n\geq k\geq0$ are integers;
\item[(b)]
$V$ is a $n$-dimensional $\mathbb C$-vector space;
\item[(c)]
$0=V_0\subset V_1\subset\cdots\subset V_{n-1}\subset V_n=V$
is a complete flag in $V$;
\item[(d)]
$x$ is an endomorphism of $V$ which leaves the flag stable and is
such that $x^2=0$;
\item[(e)]
$W$ is a $k$-dimensional subspace of $V$ such that
$\im x\subseteq W\subseteq\ker x$;
\item[(f)]
$J$ is a subset of $\{1,\ldots,n\}$ of cardinal $k$.
\end{description}

Let $\mathscr F$ be the set of all filtrations
$$\mathbf X\;:\quad
0=X_0\subseteq X_1\subseteq\cdots\subseteq X_{n-1}\subseteq X_n=W$$
such that
\begin{description}
\item[(g)]
$\dim X_p/X_{p-1}=1$ if $p\in J$ and $X_p=X_{p-1}$ if $p\notin J$.
\end{description}

We define a function $f:\mathscr F\to\mathbb Q$ as follows. Given
$\mathbf X\in\mathscr F$, we set $f(\mathbf X)=0$ except when
\begin{description}
\item[(h)]
$x(V_p)\subseteq X_p\subseteq V_p$ for any $p$.
\end{description}
When (h) holds true, we set $f(\mathbf X)=\prod_{p\in J}\eta_p$, where
$$\eta_p=\begin{cases}
1&\text{if $X_p\not\subseteq V_{p-1}$,}\\
-1&\text{if $x(V_p)\not\subseteq X_{p-1}$,}\\
0&\text{if both $X_p\subseteq V_{p-1}$ and $x(V_p)\subseteq X_{p-1}$.}
\end{cases}$$
(This definition of $\eta_p$ makes sense, because at least one of the
inclusions $X_p\subseteq V_{p-1}$ or $x(V_p)\subseteq X_{p-1}$ holds
true. In fact, if $x(V_p)\not\subseteq X_{p-1}$, then $X_p=X_{p-1}+x(V_p)$,
because $X_{p-1}$ is an hyperplane in $X_p$, and therefore
$X_p\subseteq V_{p-1}$, because both $X_{p-1}$ and $x(V_p)$ are
contained in $V_{p-1}$.)

Likewise, let $\mathscr G$ be the set of all filtrations
$$\mathbf Y\;:\quad
W=Y_0\subseteq Y_1\subseteq\cdots\subseteq Y_{n-1}\subseteq Y_n=V$$
such that
\begin{description}
\item[(j)]
$\dim Y_p/Y_{p-1}=1$ if $p\notin J$ and $Y_p=Y_{p-1}$ if $p\in J$.
\end{description}

We define a function $g:\mathscr G\to\mathbb Q$ as follows. Given
$\mathbf Y\in\mathscr G$, we set $g(\mathbf Y)=0$ except when
\begin{description}
\item[(k)]
$V_p\subseteq Y_p\subseteq x^{-1}(V_p)$ for any $p$.
\end{description}
When (k) holds true, we set $g(\mathbf Y)=\prod_{p\notin J}\eta_p$, where
$$\eta_p=\begin{cases}
1&\text{if $V_p\not\subseteq Y_{p-1}$,}\\
-1&\text{if $Y_p\not\subseteq x^{-1}(V_{p-1})$,}\\
0&\text{if both $V_p\subseteq Y_{p-1}$ and $x(Y_p)\subseteq V_{p-1}$.}
\end{cases}$$

We claim that with the above notation,
\begin{equation}
\label{eq:MainLemma}
\int_{\mathscr F}f=\int_{\mathscr G}g.
\end{equation}
The proof of this fact occupies the rest of section~\ref{se:AdapFilt}.

\subsection{}
\label{ss:AdaptBasis}
We consider the setup (a)--(e) of the previous paragraph. We show that
one can construct a basis $\mathbf e=(e_1,\ldots,e_n)$ of $V$ such
that
\begin{description}
\item[(l)]
$V_p=\linspan_{\mathbb C}\{e_1,\ldots,e_p\}$ for each $p$;
\item[(m)]
$x(e_p)\in\{0,e_1,\ldots,e_{p-1}\}$ for each $p\geq1$;
\item[(n)]
$\im x$, $\ker x$ and $W$ are coordinate subspaces w.r.t.\
the basis $\mathbf e$.
\end{description}

We will show in fact a slightly stronger statement, namely that we can
replace the subspace $W$ by an increasing sequence $\im x\subseteq W_1
\subseteq\cdots\subseteq W_\ell\subseteq\ker x$ and demand (n) for each
$W_p$ at once.

The existence of $\mathbf e$ is obvious if $V$ is a line. To prove
the general case, we use induction on $n=\dim V$.

To the sequence $W_p$, we add $W_0=\im x$ and $W_{\ell+1}=\ker x$.
Since $x$ leaves stable the hyperplane $V_{n-1}$, it induces an
endomorphism of the line $V_n/V_{n-1}$, which is necessarily zero
for $x$ is nilpotent, and therefore $V_{n-1}\supseteq\im x$. Applying
the induction hypothesis to $\widetilde V=V_{n-1}$, endowed with the
flag $(V_p)_{0\leq p\leq n-1}$, the endomorphism $\tilde x=x\bigl
|_{\widetilde V}$ and the subspaces $\widetilde W_p=W_p\cap\widetilde V$,
we get a basis $\widetilde{\mathbf e}=(e_1,\ldots,e_{n-1})$ of $V_{n-1}$
satisfying (l) and (m) for $p<n$ and such that all the $\widetilde W_p$
are coordinate subspaces w.r.t.\ the basis $\widetilde{\mathbf e}$.

We now have to complete $\widetilde{\mathbf e}$ by adding a vector
$e_n$. We distinguish two cases.

If $V_{n-1}\supseteq\ker x$, then
$x(V_{n-1})\subsetneq x(V_n)=\im x$. Now $\im x=\widetilde W_0$
is a coordinate subspace w.r.t.\ the basis $\widetilde{\mathbf e}$,
whence an index $p$ such that $e_p\in\im x\setminus x(V_{n-1})$.
We then choose $e_n\in x^{-1}(e_p)$ and check that
Conditions~(l)--(n) are satisfied.

Otherwise, there is an index $p\geq1$ such that $V_{n-1}\not\supseteq
W_p$ but $V_{n-1}\supseteq W_{p-1}$. Let us choose
$e_n\in W_p\setminus V_{n-1}$ and set $\mathbf e=(e_1,\ldots,e_n)$.
For each $q\geq p$, we have then $e_n\in W_q\setminus\widetilde W_q$,
so $W_q=\widetilde W_q+\mathbb Ce_n$, and thus $W_q$ is a coordinate
subspace w.r.t.\ the basis $\mathbf e$. Conditions~(l) and~(n) are
therefore satisfied, and Condition~(m) follows from $e_n\in\ker x$.

\subsection{}
\label{ss:CombDiag}
We consider again the setup of paragraph~\ref{ss:CombiSetup}.
Using~\ref{ss:AdaptBasis}, we find a basis $(e_1,\ldots,e_n)$
that satisfies (l)--(n). The situation can be depicted on a
diagram of the following form.
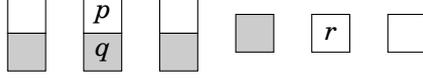
\begin{center}
\begin{tikzpicture}
\filldraw[fill=gray!40,draw=black] (0,0) rectangle (0.5,0.5);
\draw (0,0.5) rectangle (0.5,1);
\filldraw[fill=gray!40,draw=black] (1,0) rectangle (1.5,0.5);
\draw (1,0.5) rectangle (1.5,1);
\filldraw[fill=gray!40,draw=black] (2,0) rectangle (2.5,0.5);
\draw (2,0.5) rectangle (2.5,1);
\filldraw[fill=gray!40,draw=black] (3,0.25) rectangle (3.5,0.75);
\draw (4,0.25) rectangle (4.5,0.75);
\draw (5,0.25) rectangle (5.5,0.75);
\draw (1.25,0.25) node {$q$};
\draw (1.25,0.75) node {$p$};
\draw (4.25,0.5) node {$r$};
\end{tikzpicture}
\end{center}
Each column contains two indices $p,q$ such that $x(e_p)=e_q$,
an isolated box contains an index $r$ such that
$e_r\in(\ker x\setminus\im x)$, and the gray boxes are those that
contain an index $s$ such that $e_s\in W$. The above picture
represents a situation with $n=9$, $k=4$, $\dim\im x=3$,
$\dim\ker x=6$.

\subsection{}
\label{ss:CoordDesc}
To compute the left hand side of (\ref{eq:MainLemma}), we may
eliminate the filtrations $\mathbf X$ such that $f(\mathbf X)=0$,
hence we may replace $\mathscr F$ by the set $\mathscr F_0$ of
all filtrations $\mathbf X$ that satisfy (g), (h) and
\begin{description}
\item[(i)]
if $p\in J$, either $X_p\not\subseteq V_{p-1}$ or
$x(V_p)\not\subseteq X_{p-1}$.
\end{description}

We want to describe $\mathscr F_0$ in coordinates. To this aim, we
denote by $K$ the set of labels in the gray boxes; in other words,
$K=\{p\in[1,n]\mid e_p\in W\}$. In addition, let $(e_1^*,\ldots,e_n^*)$
be the dual basis to $\mathbf e$.

We first note that if $\mathscr F_0\neq\varnothing$, then
\begin{description}
\item[(A)]
Each $p\in J$ either lies in the top box of a column or belongs to $K$.
\end{description}
In fact, if $p\in J$ does not lie in the top box of a column, then
$e_p\in\ker x$. Picking $\mathbf X\in\mathscr F_0$, we thus have
$x(V_p)=x(V_{p-1})\subseteq X_{p-1}$. Condition~(i) then says that
$X_p\not\subseteq V_{p-1}$. A fortiori $W\cap V_p\not\subseteq V_{p-1}$,
and therefore $e_p\in W$.

We now claim that when Condition~(A) is satisfied, a filtration
$\mathbf X\in\mathscr F_0$ is uniquely described by a matrix
$\Xi=(\xi_{pq})$ of complex numbers with the following properties:
\begin{description}
\item[(B)]
The lines of $\Xi$ are indexed by $J$; the columns of $\Xi$ are indexed
by $K$.
\item[(C)]
If $p\in J$ does not lie in the top box of a column, then $\xi_{pp}=1$
and all the other entries in $p$-th line of $\Xi$ are zero.
\item[(D)]
If $p\in J$ lies in the top box of a column, if $q$ is the index in the
box below \raisebox{-3pt}{%
\begin{tikzpicture}
\draw (0,0) rectangle (0.5,0.5);
\draw (0.25,0.25) node {$p$};
\end{tikzpicture}%
}\,, then $\xi_{pq}=1$. In addition, if $p'\neq p$ also lies in
the top box of a column, if $q'$ is the index in the box below
\raisebox{-3pt}{%
\begin{tikzpicture}
\draw (0,0) rectangle (0.5,0.5);
\draw (0.25,0.25) node {$p'$};
\end{tikzpicture}%
}\,, and if either $p'<p$ or $p'\in J$, then $\xi_{pq'}=0$. Also, if
$q''>p$ does not lie in the top box of a column and if $q''\in J$,
then $\xi_{pq''}=0$.
\item[(E)]
Let $p\in K\setminus J$ and let $Y=(y_q)$ be the line-vector given
by $y_p=1$ and $y_q=0$ for $q\neq p$. Then $Y$ belongs to the linear
span of the rows of $\Xi$ with index $>p$.
\item[(F)]
$\Xi$ is invertible.
\end{description}

As an example, consider the following diagram
\begin{center}
\begin{tikzpicture}
\filldraw[fill=gray!40,draw=black] (0,0) rectangle (0.5,0.5);
\draw (0,0.5) rectangle (0.5,1);
\filldraw[fill=gray!40,draw=black] (1,0) rectangle (1.5,0.5);
\draw (1,0.5) rectangle (1.5,1);
\filldraw[fill=gray!40,draw=black] (2,0) rectangle (2.5,0.5);
\draw (2,0.5) rectangle (2.5,1);
\filldraw[fill=gray!40,draw=black] (3,0.25) rectangle (3.5,0.75);
\draw (4,0.25) rectangle (4.5,0.75);
\draw (5,0.25) rectangle (5.5,0.75);
\draw (0.25,0.25) node {$1$};
\draw (0.25,0.75) node {$9$};
\draw (1.25,0.25) node {$3$};
\draw (1.25,0.75) node {$4$};
\draw (2.25,0.25) node {$5$};
\draw (2.25,0.75) node {$7$};
\draw (3.25,0.5) node {$2$};
\draw (4.25,0.5) node {$6$};
\draw (5.25,0.5) node {$8$};
\end{tikzpicture}
\end{center}
and take $J=\{2,3,7,9\}$ and $K=\{1,2,3,5\}$.
Condition~(A) is satisfied.
Conditions~(B)--(D) impose $\Xi$ has the form
$$\bordermatrix{&1&2&3&5\cr
2&0&1&0&0\cr
3&0&0&1&0\cr
7&0&\xi_{72}&0&1\cr
9&1&\xi_{92}&0&0}.$$
Condition~(F) is then automatically fulfilled, while
Condition~(E) amounts to $\xi_{72}=0$.

As another example, we consider the diagram
\begin{center}
\begin{tikzpicture}
\filldraw[fill=gray!40,draw=black] (0,0) rectangle (0.5,0.5);
\draw (0,0.5) rectangle (0.5,1);
\filldraw[fill=gray!40,draw=black] (1,0) rectangle (1.5,0.5);
\draw (1,0.5) rectangle (1.5,1);
\filldraw[fill=gray!40,draw=black] (2,0) rectangle (2.5,0.5);
\draw (2,0.5) rectangle (2.5,1);
\filldraw[fill=gray!40,draw=black] (3,0) rectangle (3.5,0.5);
\draw (3,0.5) rectangle (3.5,1);
\draw (0.25,0.25) node {$1$};
\draw (0.25,0.75) node {$7$};
\draw (1.25,0.25) node {$2$};
\draw (1.25,0.75) node {$8$};
\draw (2.25,0.25) node {$3$};
\draw (2.25,0.75) node {$5$};
\draw (3.25,0.25) node {$4$};
\draw (3.25,0.75) node {$6$};
\end{tikzpicture}
\end{center}
and take $J=\{3,4,5,6\}$ and $K=\{1,2,3,4\}$.
Condition (A) is satisfied. Conditions~(B)--(D) impose $\Xi$ has the form
$$\bordermatrix{&1&2&3&4\cr
3&0&0&1&0\cr
4&0&0&0&1\cr
5&\xi_{51}&\xi_{52}&1&0\cr
6&\xi_{61}&\xi_{62}&0&1}.$$
Conditions~(E) and~(F) demand that
$\begin{vmatrix}\xi_{51}&\xi_{52}\\\xi_{61}&\xi_{62}\end{vmatrix}\neq0$.

\subsection{}
\label{ss:ProofCoordDesc}
We prove here the claim of \ref{ss:CoordDesc}.

Given the matrix $\Xi$, we set $\varphi_p=\sum_{q\in K}\xi_{pq}e_q^*$.
The correspondence $\Xi\mapsto\mathbf X$ is given by the rule
$$X_p=\{v\in W\mid\varphi_{p'}(v)=0\text{ for each $p'\in
J\cap[p+1,n]$}\}.$$

We first check that $\mathbf X$ satisfies Conditions (g)--(i)
when $\Xi$ satisfies Conditions (B)--(F).

By construction, $X_n=W$, $X_p=X_{p-1}$ if $p\notin J$, and
$X_{p-1}$ has codimension at most one in $X_p$ if $p\in J$.
In addition, $X_0=0$ thanks to (F). All this gives (g).

To establish the inclusion $X_p\subseteq V_p$, we proceed by
decreasing induction on $p$. So let us assume that $X_p\subseteq V_p$
and let us show that $X_{p-1}\subseteq V_{p-1}$. If $p\notin K$, then
$X_p\subseteq V_p\cap W\subseteq V_{p-1}$, and a fortiori
$X_{p-1}\subseteq V_{p-1}$. If $p\in J\cap K$, then $p$ does not lie
in the top box of a column, so $\varphi_p=e_p^*$ by Condition~(C),
and therefore $X_{p-1}=\{v\in X_p\mid\varphi_p(v)=0\}$ is contained
in $\{v\in V_p\mid\varphi_p(v)=0\}=V_{p-1}$. Lastly, if
$p\in K\setminus J$, then $e_p^*\bigl|_W$ is a linear combination of
the forms $\varphi_{p'}\bigl|_W$, for $p'\in J\cap[p+1,n]$, thanks to
Condition~(E); it follows that $e_p^*$ vanishes on $X_p$, hence that
$X_p$ is already contained in $V_{p-1}$.

Let now $p$ be in the top box of a column and let $q$ be the index in
the box below \raisebox{-3pt}{%
\begin{tikzpicture}
\draw (0,0) rectangle (0.5,0.5);
\draw (0.25,0.25) node {$p$};
\end{tikzpicture}%
}\,. Take $p'\in J\cap[p+1,n]$. If $p'$ does not lie in the top box
of a column, then $\varphi_{p'}=e_{p'}^*$ and $p'>p>q$. If $p'$
lies in the top box of a column, then $\xi_{p'q}=0$ by Condition~(D).
In either case, $\varphi_{p'}(e_q)=0$. Since this holds for any
$p'\in J\cap[p+1,n]$, we conclude that $e_q\in X_p$. This can be
rewritten $x(e_p)\in X_p$. In fact, $x(e_p)\in X_p$ holds for any
$p$, because $x(e_p)=0$ when $p$ is not in a top box. It follows that
$x(V_p)\subseteq X_p$. Combined with the inclusion $X_p\subseteq V_p$
proved in the previous paragraph, we get~(h).

Fix $p\in J$. If $p$ does not lie in the top box of a column,
then $X_{p-1}=\{v\in X_p\mid e_p^*(v)=0\}$ by Condition~(C);
in addition, we have already proved (g), so we know that
$X_{p-1}\neq X_p$; these two facts imply that $X_p$ cannot be
contained in $\ker e_p^*$, that is, $X_p\not\subseteq V_{p-1}$.
If $p$ lies in the top box of a column, then $\varphi_p(x(e_p))=1$
by Condition~(D), so $x(e_p)\notin X_{p-1}$, and therefore
$x(V_p)\not\subseteq X_{p-1}$. This shows (i).

Therefore the map $\Xi\mapsto\mathbf X$ is well-defined. To show its
bijectivity, we construct its inverse. We thus start with a
filtration $\mathbf X\in\mathscr F_0$ and look for the matrix $\Xi$.

Take $p\in J$. The $p$-th line of $\Xi$ displays the coordinates of
a linear form $\varphi_p\in\linspan_{\mathbb C}\{e_q^*\mid q\in K\}$
such that $X_{p-1}=\{v\in X_p\mid\varphi_p(v)=0\}$. We need to
show the existence of $\varphi_p$ and to normalize it in a unique
way. By induction, we may assume that the $p'$-th line of $\Xi$ has
been determined for each $p'\in J$ greater than $p$.

Assume first that $p$ does not lie in the top box of a column. Then
by (h), we have $x(V_p)=x(V_{p-1})\subseteq X_{p-1}$, and thus
$X_p\not\subseteq V_{p-1}$, thanks to Condition (i). It follows that
$V_{p-1}\cap X_p$ is strictly contained in $X_p$. However $X_{p-1}$
is contained in $V_{p-1}\cap X_p$ and is an hyperplane of $X_p$, by
(g), so we necessarily have $X_{p-1}=V_{p-1}\cap X_p$. Therefore
$X_{p-1}$ is obtained by cutting $X_p$ by the hyperplane of equation
$\varphi_p=e_p^*$. This gives Condition (C) on the matrix $\Xi$.

Now assume that $p$ lies in the top box of a column. Then
$V_{p-1}\cap W=V_p\cap W$, so $V_{p-1}\cap X_p=V_p\cap X_p=X_p$,
whence $X_p\subseteq V_{p-1}$. Condition~(i) then says that
$x(V_p)\not\subseteq X_{p-1}$. Looking at Conditions (g) and (h),
we deduce that $X_{p-1}$ is produced by cutting $X_p$ by an
hyperplane that contains $x(V_{p-1})$ but not $x(V_p)$. If $q$
is the index in the box below \raisebox{-3pt}{%
\begin{tikzpicture}
\draw (0,0) rectangle (0.5,0.5);
\draw (0.25,0.25) node {$p$};
\end{tikzpicture}%
}\,, then the equation of this hyperplane can be uniquely written as
$\varphi_p=e_q^*+\sum_{r\neq q}\xi_{pr}e_r^*$. If $p'<p$ lies in the
top box of a column and if $q'$ is the index of the box below
\raisebox{-3pt}{%
\begin{tikzpicture}
\draw (0,0) rectangle (0.5,0.5);
\draw (0.25,0.25) node {$p'$};
\end{tikzpicture}%
}\,, then $\xi_{pq'}=0$, because
$x(e_{p'})\in x(V_{p-1})$ and $x(V_{p-1})\subseteq\ker\varphi_p$. We
now modify $\varphi_p$ so as to obtain the form specified in Condition (D),
without altering the property $X_{p-1}=\{v\in X_p\mid\varphi_p(v)=0\}$:
\begin{itemize}
\item
Since $X_p\subseteq W$, we may subtract from $\varphi_p$ all the
terms $\xi_{pr}e_r^*$ such that $r\notin K$.
\item
If $p'>p$ lies in the top box of a column and if $p'\in J$, then
$X_p\subseteq X_{p'-1}\subseteq\ker\varphi_{p'}$. Let $q'$ be the
index in the box below \raisebox{-3pt}{%
\begin{tikzpicture}
\draw (0,0) rectangle (0.5,0.5);
\draw (0.25,0.25) node {$p'$};
\end{tikzpicture}%
}\,; we may then subtract $\xi_{pq'}\varphi_{p'}$ from $\varphi_p$,
which yields a new $\varphi_p$ with $\xi_{pq'}=0$.
\item
We repeat the previous step for each relevant $p'$ in order to
annihilate all $\xi_{pq'}$.
\item
Lastly, if $q''>p$ does not lie a the top box of a column and if
$q''\in J$, then $X_p\subseteq V_p\subseteq\ker e_{q''}^*$, so we
may subtract the term $\xi_{pq''}e_{q''}^*$ from $\varphi_p$.
\end{itemize}

It remains to check that the matrix $\Xi$ satisfies (E) and (F).
We note that by construction,
$$X_p=\{v\in W\mid\varphi_{p'}(v)=0\text{ for each $p'\in
J\cap[p+1,n]$}\}.$$
Condition (F) then comes from $X_0=0$. Now take $p\in K\setminus J$.
By (g) and (h), we have $X_p=X_{p-1}\subseteq V_{p-1}$, so $e_p^*$
vanishes on $X_p$. From the above description of $X_p$, it then
follows that $e_p^*\bigl|_W$ belongs to the linear span of the elements
$\varphi_{p'}\bigl|_W$ with $p'\in J\cap[p+1,n]$. This is Condition (E).

\subsection{}
\label{ss:EulerCharF0}
Sections~\ref{ss:CoordDesc} and~\ref{ss:ProofCoordDesc} gives a
description of $\mathscr F_0$ as a set of matrices. This allows to
compute the Euler characteristic of $\mathscr F_0$: it is equal to
zero or one, the latter case happening precisely when each column
of the diagram contains exactly one element of $J$ and when the
remaining elements of $J$ occupy the isolated gray boxes. The aim
of this section is to show this result.

We begin with a closer look at the Conditions (B)--(F). Conditions
(C) and (D) prescribe that certain matrix elements $\xi_{pq}$ have
value zero or one. The other matrix entries of $\Xi$ can be freely
chosen; we call them the free entries. The set $\mathscr V$ of
all matrices defined by Conditions (B)--(D) is thus a finite
dimensional $\mathbb C$-vector space, spanned by the free entries.
Conditions (E) and (F) define a locally closed subset
$\mathscr G\subseteq\mathscr V$.

The position of the zeros and the ones in a matrix $\Xi\in\mathscr V$
obey a special pattern, which can be described as follows. Up to
a reordering of the rows, $\Xi$ can be viewed as the pile of two matrices
$\Xi_1$ and $\Xi_2$. The matrix $\Xi_1$ gathers the rows with index in
$J\cap K$; by Condition~(C), each row of $\Xi_1$ is a basis row-vector,
that is, its entries are all zero except one entry equal to one. The
matrix $\Xi_2$ gather the rows with index in $J\setminus K$; all free
entries of $\Xi$ are in $\Xi_2$. Condition~(D) implies that a column of
$\Xi_2$ either is a basis column-vector, or its entries are zero or free
entries. In other words, if a column of $\Xi_2$ contains a one, then it
is a basis column-vector.

With the help of this description, one easily sees that a minor of
$\Xi$ is a homogeneous polynomial in the free entries. In fact, let
$H$ be a matrix obtained from $\Xi$ by removing some lines and some
columns. Like $\Xi$, the matrix $H$ has a description as a pile of
two matrices $H_1$ and $H_2$; the only difference with $\Xi$ is that
a row of $H_1$ is either a basis row-vector or zero. To compute the
determinant of $H$, we pick a row in $H_1$. If the row is zero, then
the determinant vanishes; otherwise, this row contains a one, and
the determinant is changed at most by a sign if we remove the line and
the column that contains this one. Repeting this operation, we get rid
of all the lines of $H_1$ and end up with a matrix obtained by removing
columns from $H_2$. Again, a column of this matrix either is a basis
column-vector, or its entries are zero or free entries. Pursuing
the expansion of the determinant, we remove the basis vector-columns
together with some lines. We wind up with a matrix $H_3$ whose entries
are zeros or free entries, and we have $\det H=\pm\det H_3$. This
establishes our claim.

In view of this property, Conditions (E) and (F) define
$\mathscr G$ as a cone in $\mathscr V$. Therefore
$\mathscr G_0=\mathscr G\setminus\{0\}$ is endowed with a
free action of $\mathbb C^*$. The principal $\mathbb C^*$-bundle
$\mathscr G_0\to\mathscr G_0/\mathbb C^*$ then gives for the Euler
characteristic with compact support
$$\chi(\mathscr G_0)=\chi(\mathscr G_0/\mathbb C^*)\;\chi(\mathbb C^*)
=0.$$
We conclude that $\chi(\mathscr F_0)=\chi(\mathscr G)$ is zero
if $\mathscr G=\mathscr G_0$ and is one if
$\mathscr G=\mathscr G_0\sqcup\{0\}$. To finish the proof, it remains
to determine when the origin of the vector space $\mathscr V$ belongs
to $\mathscr G$.

According to the computation explained above, $\det\Xi=\pm\det\Xi_3$,
where the entries of $\Xi_3$ are either zero or free. In view of
Condition~(F), the origin of $\mathscr V$ belongs to $\mathscr G$ if
and only if $\Xi_3$ is the empty matrix. A scrutiny of the construction
of $\Xi_3$ reveals that the columns of $\Xi_3$ are indexed by columns
\raisebox{-11pt}{%
\begin{tikzpicture}
\filldraw[fill=gray!40,draw=black] (0,0) rectangle (0.5,0.5);
\draw (0,0.5) rectangle (0.5,1);
\draw (0.25,0.25) node {$q$};
\draw (0.25,0.75) node {$p$};
\end{tikzpicture}%
} with $p,q\notin J$ and boxes \raisebox{-3pt}{%
\begin{tikzpicture}
\filldraw[fill=gray!40,draw=black] (0,0) rectangle (0.5,0.5);
\draw (0.25,0.25) node {$r$};
\end{tikzpicture}%
} with $r\notin J$. Therefore $\Xi_3$ is the empty matrix if and only if
each column of the diagram contains at least one element of $J$ and
each isolated gray box contains an element of $J$. By cardinality, the
latter condition is equivalent to the criterion given in our claim.

\subsection{}
\label{ss:ProofMainLemma}
We are now in a position to prove (\ref{eq:MainLemma}).

We first consider the left hand side of the equality. By definition
of $f$, given $\mathbf X\in\mathscr F$, we have $f(\mathbf X)\neq0$
if and only if $\mathbf X\in\mathscr F_0$. In addition, in the case
$\mathbf X\in\mathscr F_0$, we have $x(V_p)\not\subseteq X_{p-1}$
if and only if $p$ lies in the top box of the diagram,
by section~\ref{ss:ProofCoordDesc}. It follows that $f$ assumes only
two values on $\mathscr F$: it vanishes on
$\mathscr F\setminus\mathscr F_0$ and it is equal to $(-1)^N$ on
$\mathscr F_0$, where $N$ is the number of indices $p\in J$ in the
top boxes of columns. Therefore
$\int_{\mathscr F}f=(-1)^N\chi(\mathscr F_0)$.
The computation in section~\ref{ss:EulerCharF0} then gives us the
following combinatorial recipe for the value of $\int_{\mathscr F}f$:
\begin{itemize}
\item
It is zero if there is an isolated white box whose index belongs to
$J$ or if there is a column whose two indices belong to $J$;
\item
Otherwise, it is $(-1)^N$, where $N=\bigl|J\setminus K\bigr|$.
\end{itemize}

The right hand side $\int_{\mathscr G}g$ can be computed by duality.
Indeed, we define a twisted setup as follows:
\begin{description}
\item[(\~a)]
$\tilde n=n$, $\tilde k=n-k$;
\item[(\~b)]
$\widetilde V$ is the dual of $V$;
\item[(\~c)]
The $p$-dimensional subspace $\widetilde V_p$ of the flag of
$\widetilde V$ is the orthogonal of $V_{n-p}$;
\item[(\~d)]
$\tilde x$ is the transpose of $x$;
\item[(\~e)]
$\widetilde W$ is the orthogonal of $W$;
\item[(\~f)]
$\widetilde J=\{n+1-p\mid p\in[1,n]\setminus J\}$.
\end{description}
In addition, duality also gives us a bijective correspondence between
$\mathbf Y$ and the set $\widetilde{\mathbf X}$ of filtrations for the
twisted setup; the map $g$ is thereby transported to a map $\tilde f$.
We then have $\int_{\mathscr G}g=\int_{\widetilde{\mathscr F}}\tilde f$.

We can then use again our combinatorial recipe. We observe that the
diagram of the twisted situation is obtained by turning upside-down
the diagram in section~\ref{ss:CombDiag}, by changing each entry $p$
to $n+1-p$, and by inverting the colors in all the boxes. Comparing
this with the change described in (\~f), we retrieve for
$\int_{\widetilde{\mathscr F}}\tilde f$ the same value as for
$\int_{\mathscr F}f$. Therefore both sides of~(\ref{eq:MainLemma})
are equal.

\section{The preprojective model for $\mathbf U^+$}
\label{se:PrepModel}
\subsection{}
\label{ss:EnvAlg}
We consider a finite non-empty graph without loops. This is the same as
giving two finite sets $I$ and $H$ with $I\neq\varnothing$, a fixed point
free involution $h\mapsto h^*$ of $H$, and two maps $s,t:H\to I$ such
that $s(h^*)=t(h)\neq s(h)$ for all $h\in H$.

For $i,j\in I$, we set $a_{ij}=-\bigl|\{h\in H\mid s(h)=i,\;t(h)=j\}\bigr|$
if $i\neq j$ and $a_{ij}=2$ if $i=j$. Then $A=(a_{ij})$ is a symmetric
generalized Cartan matrix.

Let $\mathbf U^+$ be the $\mathbb Q$-algebra defined by generators $e_i$
for $i\in I$, submitted to the Serre relations
$$\sum_{p+q=1-a_{ij}}(-1)^p\frac{e_i^p}{p!}\;e_j\;\frac{e_i^q}{q!}=0.$$

The weight lattice is the free $\mathbb Z$-module with basis
$\{\alpha_i\mid i\in I\}$; we denote it by $\mathbb ZI$. We write
a typical element in $\mathbb ZI$ as $\nu=\sum_{i\in I}\nu_i\alpha_i$;
when all the $\nu_i$ are non-negative, we write $\nu\in\mathbb NI$.
Given $\nu\in\mathbb NI$, we denote by $\mathbf U^+_\nu$ the subspace
of $\mathbf U^+$ spanned by the monomials $e_{i_1}\cdots e_{i_n}$ for
sequences $i_1,\ldots,i_n$ in which $i$ appears $\nu_i$ times.

\subsection{}
\label{ss:SubalgUi}
For $i\neq j$ and $0\leq m\leq-a_{ij}$, let
$$f_{i,j,m}=\sum_{p+q=m}(-1)^p\frac{e_i^p}{p!}\;e_j\;\frac{e_i^q}{q!}.$$
For a fixed $i\in I$, we denote by $\mathbf U^+_i$ the subalgebra of
$\mathbf U^+$ generated by the elements $f_{i,j,m}$, for all possible
$(j,m)$. We set $\mathbf U^+_{i,\nu}=\mathbf U^+_i\cap\mathbf U^+_\nu$.

We define an automorphism $T_i$ of $\mathbf U^+_i$ by the requirement
$$T_i(f_{i,j,m})=(-1)^mf_{i,j,-a_{ij}-m}.$$
Lemma~38.1.3 in \cite{Lusztig93} shows that this definition makes sense;
in fact, $T_i$ coincides with the restriction to $\mathbf U^+_i$ of
(the specialization at $v=1$ of) the automorphism $T'_{i,-1}$.
Lemma~38.1.2 in \cite{Lusztig93} then implies that the subalgebra
$\mathbf U^+_i$ just defined coincides with the subalgebra $\mathbf U^+_i$
of the introduction.

It is here worth noting that $\mathbf U^+_i$ is the smallest subalgebra
of $\mathbf U^+$ that contains all the elements $e_j$ for $j\neq i$ and
that is stable by the derivation $D_i=\ad(e_i)$. In fact, an immediate
computation based on the binomial theorem shows that
$f_{i,j,m}=\frac{(-1)^m}{m!}D^m(e_j)$.

\subsection{}
\label{ss:NilVar}
We fix a function $\varepsilon:H\to\{\pm1\}$ such that
$\varepsilon(h)+\varepsilon(h^*)=0$ for all $h\in H$.
We view the tuple $(I,H,s,t)$ as directed graph, $s$ and $t$ being
the source and target maps. We work over the field of complex numbers.

The preprojective algebra is defined as the quotient of the path
algebra of this graph by the ideal generated by
$$\sum_{h\in H}\varepsilon(h)h^*h.$$
Thus, a finite dimensional representation $M$ of the preprojective
algebra is the datum of finite dimensional $\mathbb C$-vector
spaces $M_i$ together with a tuple
$$(M_h)\in\prod_{h\in H}\Hom_{\mathbb C}(M_{s(h)},M_{t(h)})$$
such that for any $i\in I$,
$$\sum_{\substack{h\in H\\s(h)=i}}\varepsilon(h)M_{h^*}M_h=0$$
as a linear map $M_i\to M_i$. The dimension-vector of $M$ is defined as
$$\dimvec M=\sum_{i\in I}(\dim M_i)\alpha_i.$$

Each vertex $i\in I$ affords a one-dimensional representation $S_i$
of the preprojective algebra; explicitly, $\dimvec S_i=\alpha_i$
and the arrows act by zero on $S_i$. A finite dimensional
representation $M$ of the preprojective algebra is said to be
nilpotent if all its Jordan-H\"older components are one-dimensional,
that is, are isomorphic to a $S_i$. Nilpotent representations of
the preprojective algebra are the same as finite dimensional
$\Pi$-modules, where $\Pi$ is the completion of the preprojective
algebra with respect to the ideal generated by the arrows.

We can here fix the vector space $M_i$ and let the linear maps $M_h$
vary: we denote by $\mathscr C$ the category of finite dimensional
$I$-graded $\mathbb C$-vector spaces $\mathbf V=\bigoplus_{i\in I}V_i$,
and for given $\mathbf V\in\mathscr C$, we denote by $\Lambda_{\mathbf V}$
the set of all elements $x=(x_h)$ in
$$E_{\mathbf V}=\prod_{h\in H}\Hom_{\mathbb C}(V_{s(h)},V_{t(h)})$$
such that $(\mathbf V,x)$ is a nilpotent representation of the
preprojective algebra. An isomorphism $\mathbf V\cong\mathbf W$
in $\mathscr C$ induces an isomorphism
$\Lambda_{\mathbf V}\cong\Lambda_{\mathbf W}$; in particular, the
group $G_{\mathbf V}=\prod_{i\in I}\GL(V_i)$ acts on
$\Lambda_{\mathbf V}$. Isomorphism classes of $\Pi$-modules of
dimension-vector $\dimvec\mathbf V$ are in one-to-one correspondence
with $G_{\mathbf V}$-orbits in~$\Lambda_{\mathbf V}$.

\subsection{}
\label{ss:ConstFunNilVar}
Recall the notation defined in section~\ref{ss:ConstFun}.

For $\mathbf V\in\mathscr C$, let $\widetilde M(\Lambda_{\mathbf V})$
be the space of all functions $f\in M(\Lambda_\mathbf V)$ that are
constant on any $G_{\mathbf V}$-orbit in $\Lambda_{\mathbf V}$. If
$\mathbf V,\mathbf W\in\mathscr C$ have the same dimension-vector,
say $\nu$, then $\widetilde M(\Lambda_{\mathbf V})$ and
$\widetilde M(\Lambda_{\mathbf W})$ are canonically isomorphic. One
can therefore identify these spaces and safely denote them by
$\widetilde M_\nu$.

Let $\nu',\nu''\in\mathbb NI$ and set $\nu=\nu'+\nu''$. We define a
bilinear map $\star:\widetilde M_{\nu'}\times\widetilde M_{\nu''}\to
\widetilde M_\nu$ by the following recipe. We choose
$\mathbf V,\mathbf V',\mathbf V''\in\mathscr C$ such that
$\nu=\dimvec\mathbf V$, $\nu'=\dimvec\mathbf V'$ and
$\nu''=\dimvec\mathbf V''$. For $(f',f'')\in\widetilde
M(\Lambda_{\mathbf V'})\times\widetilde M(\Lambda_{\mathbf V''})$
and $x\in\Lambda_{\mathbf V}$, we set $(f'\star f'')(x)=\int_{\mathscr
H}\phi$, where the following notation is used:
\begin{itemize}
\item
$\mathscr H$ is the variety consisting of all $I$-graded
subspaces $\mathbf W\subseteq\mathbf V$ such that
$\dimvec\mathbf W=\dimvec\mathbf V''$ (this is a product
of Grassmannian varieties).
\item
Given $\mathbf W\in\mathscr H$, the number $\phi(\mathbf W)$ is zero
except if $x_h(W_{s(h)})\subseteq W_{t(h)}$ for all $h\in H$.
In the latter case, let $\tilde x'\in\Lambda_{\mathbf V/\mathbf W}$
and $\tilde x''\in\Lambda_{\mathbf W}$ be the elements induced by $x$,
and let $x'\in\Lambda_{\mathbf V'}$ and $x''\in\Lambda_{\mathbf V''}$
be the elements obtained by transporting $\tilde x'$ and $\tilde x''$
through isomorphisms $\mathbf V'\cong\mathbf V/\mathbf W$ and
$\mathbf V''\cong\mathbf W$; then $\phi(\mathbf W)=f'(x')f''(x'')$.
\end{itemize}

The maps $\star$ combine to endow the $\mathbb Q$-vector space
$\widetilde M=\bigoplus_{\nu\in\mathbb NI}\widetilde M_\nu$ with
the structure of a $\mathbb NI$-graded algebra.

With this notation, there is a unique morphism of algebras
$\kappa:\mathbf U^+\to\widetilde M$ such that for any $i\in I$ and
any $p\geq 1$, for any $\mathbf V\in\mathscr C$ of dimension-vector
$p\alpha_i$, the element $\kappa(e_i^p/p!)$ is the function
$f\in\widetilde M(\Lambda_{\mathbf V})$ with value $1$ on the
point $\Lambda_{\mathbf V}$. The morphism $\kappa$ is injective
(Theorem~2.7~(c) in \cite{Lusztig00}).

\subsection{}
\label{ss:ReflFunc}
In section~2.2 of \cite{BaumannKamnitzer10}, endofunctors $\Sigma_i$
and $\Sigma_i^*$ of the category $\Pi\mmod$, called reflection
functors, are defined for each vertex $i\in I$. In this paper, we
will use only $\Sigma_i^*$; it can be quickly defined as
$\Sigma_i^*=I_i\otimes_\Pi?$, where $I_i$ is the annihilator of the
simple $\Pi$-module $S_i$.

To concretely describe $\Sigma_i^*$, we introduce a special notation,
that analyzes a $\Lambda$-module $M$ locally around the vertex $i$. We
break the datum of $M$ in two parts: the first part
consists of the vector spaces $M_j$ for $j\neq i$ and of the linear
maps between them; the second part consists of the vector spaces and
of the linear maps that appear in the diagram
$$\bigoplus_{\substack{h\in H\\s(h)=i}}M_{t(h)}
\xrightarrow{\ (M_{h^*})\ }M_i\xrightarrow{(\varepsilon(h)M_h)}
\bigoplus_{\substack{h\in H\\s(h)=i}}M_{t(h)}.$$
For brevity, we will write the latter as
\begin{equation}
\label{eq:LocPresMod}
\widetilde M_i\xrightarrow{M_{\iin(i)}}M_i
\xrightarrow{M_{\out(i)}}\widetilde M_i.
\end{equation}
The relations of the preprojective algebra imply
\begin{equation}
\label{eq:MinMout}
M_{\iin(i)}M_{\out(i)}=0.
\end{equation}

With this notation, the module $\Sigma_i^*M$ is obtained by
replacing~(\ref{eq:LocPresMod}) with
$$\widetilde M_i\twoheadrightarrow\coker M_{\out(i)}
\xrightarrow{M_{\out(i)}\overline M_{\iin(i)}}\widetilde M_i,$$
where the map $\overline M_{\iin(i)}:\coker M_{\out(i)}\to M_i$ is
induced by $M_{\iin(i)}$. The vector spaces $M_j$ for $j\neq i$ and
the linear maps between them are not affected by the functor $\Sigma_i^*$.

\subsection{}
\label{ss:GenForm}
Let $\nu\in\mathbb NI$ and let $\mathbf V\in\mathscr C$ be of
dimension-vector $\nu$. The evaluation at a point
$x\in\Lambda_{\mathbf V}$ is a linear form on
$\widetilde M(\Lambda_{\mathbf V})=\widetilde M_\nu$, hence gives
via $\kappa$ a linear form $\delta_x$ on $\mathbf U^+_\nu$. If $M$
is a $\Pi$-module of dimension-vector $\nu$, we write $\delta_M$
instead of $\delta_x$, where $x\in\Lambda_{\mathbf V}$ is chosen
so that $(\mathbf V,x)$ is isomorphic to $M$.

In section~\ref{se:ProofGenForm}, we will show the following
statement. Let $M$ be a $\Pi$-module of dimension-vector~$\nu$.
If $\ker M_{\out(i)}=0$, then
\begin{equation}
\label{eq:GenForm}
\forall u\in\mathbf U^+_{i,\nu},\qquad\langle\delta_M,u\rangle
=\langle\delta_{\Sigma_i^*M},T_i(u)\rangle.
\end{equation}

In the remainder of section~\ref{se:PrepModel}, we explain how to
deduce (\ref{eq:DualForm}) from this result.

\subsection{}
\label{ss:KSCrystal}
Given an algebraic variety $X$, we denote by $\Irr X$ the set of
irreducible components of~$V$.

Let $\nu\in\mathbb NI$. Given $\mathbf V,\mathbf W\in\mathscr C$
of dimension-vector $\nu$, one can construct an isomorphism
$\Lambda_{\mathbf V}\cong\Lambda_{\mathbf W}$, determined
up to composition with an element of $G_{\mathbf V}$. The latter
being connected, one gets a canonical bijection
$\Irr\Lambda_{\mathbf V}\cong\Irr\Lambda_{\mathbf W}$. One can
therefore identify these sets and safely denote them by
$B_\nu$.

In section~8 of \cite{Lusztig90}, Lusztig endows the set
$B=\bigsqcup_{\nu\in\mathbb NI}B_\nu$ with the structure of a crystal
in the sense of Kashiwara (see section~3 in \cite{KashiwaraSaito97}
for a summary review of this notion). Through the study of the
properties of an involution $b\mapsto b^*$, Kashiwara and Saito prove
(Theorem~5.3.2 in \cite{KashiwaraSaito97}) that $B$ is isomorphic to
the crystal $B(-\infty)$ associated to the crystal basis of
$U_q^+(\mathfrak g)$. In this model, the values $\varphi_i(b)$
and $\varphi_i(b^*)$ are given as follows: if $b\in B(-\infty)$
corresponds to the irreducible component $Z\in\Irr\Lambda_{\mathbf V}$
under Kashiwara and Saito's isomorphism and if $x$ is general enough
in $Z$, then
$$\varphi_i(b)=\dim\coker M_{\iin(i)}\quad\text{and}\quad
\varphi_i(b^*)=\dim\ker M_{\out(i)},$$
where $M=(\mathbf V,x)$.

In addition, Saito's reflection $\sigma_i$ (which we define in
section~\ref{ss:Corollaries}) also has a nice interpretation:
if $b\in B_i$ corresponds to $Z'\in\Irr\Lambda_{\mathbf V'}$,
if $\sigma_ib$ corresponds to $Z''\in\Irr\Lambda_{\mathbf V''}$,
then any non-empty open subset of $Z'\times Z''$ contains an element
$(x',x'')$ such that
$$\Sigma_i^*(\mathbf V',x')\cong(\mathbf V'',x'').$$
For a proof, see Proposition~18 in~\cite{BaumannKamnitzer10}.

\subsection{}
\label{ss:SemicanBas}
Let $\nu\in\mathbb NI$, let $\mathbf V\in\mathscr C$ be of
dimension-vector $\nu$, and let $Z\in\Irr\Lambda_{\mathbf V}$.
Given $u\in\mathbf U^+_\nu$, the function
$\kappa(u)=(x\mapsto\delta_x(u))$ is constructible on
$\Lambda_{\mathbf V}$, so takes a constant value on a non-empty
open subset of $Z$. We denote this value by $\delta_Z(u)$.
Since $\mathbf U^+_\nu$ is finite dimensional, the open subset
can be chosen independently of $u$: there is a non-empty open
subset $\Omega\subseteq Z$ such that for any $x\in\Omega$, we
have~$\delta_x=\delta_Z$ as linear forms
$\mathbf U^+_\nu\to\mathbb Q$.

By Theorem~2.7 in~\cite{Lusztig00},
$\{\delta_Z\mid Z\in\Irr\Lambda_{\mathbf V}\}$ is
a basis of $(\mathbf U^+_\nu)^*$. This basis does not depend on
the choice of $\mathbf V$ and is called the dual semicanonical
basis. If $Z$ corresponds to $b\in B_\nu$ under Kashiwara and
Saito's bijection, then the element $S(b)^*$ used in
section~\ref{ss:Corollaries} is $S(b)^*=\delta_Z$.

Now take $b\in B_i$. Identify $b$ to an irreducible component
$Z'\in\Irr\Lambda_{\mathbf V'}$ and identify $\sigma_ib$ to an
irreducible component $Z''\in\Irr\Lambda_{\mathbf V''}$, as we
did in section~\ref{ss:KSCrystal}. Let $\Omega'\subseteq Z'$ be a
non-empty open subset such that $\delta_{x'}=\delta_{Z'}$ for any
$x'\in\Omega'$. Likewise, let $\Omega''\subseteq Z''$ be a non-empty
open subset such that $\delta_{x''}=\delta_{Z''}$ for any
$x''\in\Omega''$. By shrinking $\Omega'$ if necessary, we can
assume that
$$0=\varphi_i(b^*)=\dim\ker M_{\out(i)}$$
for any $x'\in\Omega'$, where $M=(\mathbf V',x')$.

Take $(x',x'')\in\Omega'\times\Omega''$ such that
$$\Sigma_i^*(\mathbf V',x')\cong(\mathbf V'',x'').$$
Applying Equation (\ref{eq:GenForm}) to $M=(\mathbf V',x')$, we get
$$\forall u\in\mathbf U^+_{i,\nu},\qquad\langle\delta_{Z'},u\rangle
=\langle\delta_{Z''},T_i(u)\rangle,$$
where $\nu$ is the weight of $b$. This equation is
(\ref{eq:DualForm}), with $b$ instead of $b'$.
In other words, we have showed that (\ref{eq:DualForm}) is a
corollary to (\ref{eq:GenForm}).

\section{Proof of (\ref{eq:GenForm})}
\label{se:ProofGenForm}
It remains to show (\ref{eq:GenForm}). This is the purpose of this
section.

\subsection{}
\label{ss:StarDiag}
We first look at a particular case of (\ref{eq:GenForm}). Namely, we
consider the following star-shaped graph with set of vertices
$I=\{0,\ldots,n\}$.
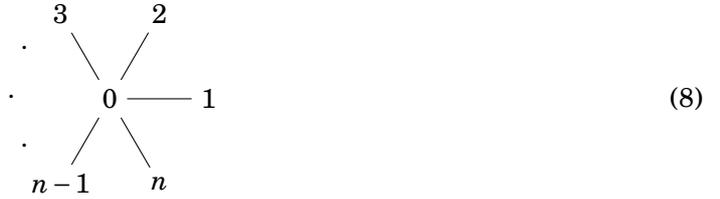
\begin{equation}
\label{eq:StarDiag}
\raisebox{-36pt}{%
\begin{tikzpicture}
\node (0) at (0,0) {$0$};
\node (1) at (0:1.3) {$1$};
\node (2) at (60:1.3) {$2$};
\node (3) at (120:1.3) {$3$};
\node at (150:1.3) {$\cdot$};
\node at (180:1.3) {$\cdot$};
\node at (210:1.3) {$\cdot$};
\node (nmu) at (240:1.3) {$n-1$};
\node (n) at (300:1.3) {$n$};
\draw (0) -- (1);
\draw (0) -- (2);
\draw (0) -- (3);
\draw (0) -- (nmu);
\draw (0) -- (n);
\end{tikzpicture}%
}
\end{equation}

Let $0\leq k\leq n$ and let $M$ be a $\Pi$-module of dimension-vector
$\nu=k\alpha_0+(\alpha_1+\cdots+\alpha_n)$. We claim that if
$\ker M_{\out(0)}=0$, then (\ref{eq:GenForm}) holds for $i=0$.

By linearity, it suffices to check (\ref{eq:GenForm}) for $u$ of the
form $u=f_{0,\pi(n),m_n}f_{0,\pi(n-1),m_{n-1}}\cdots f_{0,\pi(1),m_1}$,
where $\pi$ is a permutation of $[1,n]$ and where each $m_j\in\{0,1\}$.
The condition that $u$ has weight $\nu$ imposes that
$J=\{j\in[1,n]\mid m_j=1\}$ has cardinal $k$.

Let us set
$$V=\widetilde M_0=M_1\oplus\cdots\oplus M_n,\quad
W=\im M_{\out(0)}\quad\text{and}\quad x=M_{\out(0)}M_{\iin(0)}.$$
The preprojective relation $M_{\iin(0)}M_{\out(0)}=0$ ensures that
$x^2=0$ and that $W\subseteq\ker x$. We further set $V_0=0$ and
$V_p=M_{\pi(1)}\oplus\cdots\oplus M_{\pi(p)}$ for each $p\in[1,n]$;
this gives a complete flag in $V$.

If $x$ does not leaves this flag stable, then both sides of
(\ref{eq:GenForm}) vanish, so (\ref{eq:GenForm}) holds true.

If $x$ leaves this flag stable, then we are in the setup of
section~\ref{ss:CombiSetup}. Then the two sides of (\ref{eq:GenForm})
are given by the integrals on the two sides of (\ref{eq:MainLemma}).
The proof of~(\ref{eq:MainLemma}) in section~\ref{se:AdapFilt}
therefore ensures that our particular case of (\ref{eq:GenForm}) holds true.

\subsection{}
\label{ss:PartInt}
We now tackle the general case of (\ref{eq:GenForm}).

Since $\mathbf U^+_i$ is generated by the elements $f_{i,j,m}$,
it is enough to consider the case of a monomial
$u=f_{i,j_r,m_r}\cdots f_{i,j_1,m_1}$, where $j_s\in I\setminus\{i\}$
and $0\leq m_s\leq-a_{ij_s}$ for each $s\in[1,r]$. The weight of $u$
is of course $\nu=m\alpha_i+\alpha_{j_1}+\cdots+\alpha_{j_r}$, where
$m=m_1+\cdots+m_r$. We consider a $\Pi$-module $M$ of dimension-vector
$\nu$ such that $\ker M_{\out(i)}=0$. We write $M=(\mathbf V',x')$
and $\Sigma_i^*M=(\mathbf V'',x'')$, where $\mathbf V'\in\mathscr C$,
$\mathbf V''\in\mathscr C$, $x'\in\Lambda_{\mathbf V'}$ and
$x''\in\Lambda_{\mathbf V''}$. In fact, the reflection functor
$\Sigma_i^*$ only touches the space attached to vertex $i$, so for
$k\neq i$, we may identify $V''_k$ to $V'_k$ and write simply $V_k$.

The left hand side of (\ref{eq:GenForm}) is the evaluation at $M$ of the
function $\kappa(f_{i,j_r,m_r})\star\cdots\star\kappa(f_{i,j_1,m_1})$.
By the definition in section~\ref{ss:ConstFunNilVar}, this number is
an integral $\int_{\mathscr H'}\phi'$, where $\mathscr H'$ is the set
of all filtrations
$$0=\mathbf V'_0\subseteq\mathbf V'_1\subseteq\cdots\subseteq\mathbf
V'_{r-1}\subseteq\mathbf V'_r=\mathbf V'$$
such that $\dimvec(\mathbf V'_s/\mathbf V'_{s-1})=\alpha_{j_s}+m_s\alpha_i$
for each $s\in[1,r]$ and where $\phi':\mathscr H\to\mathbb Q$ is a
function given by the product of the~$\kappa(f_{i,j_s,m_s})$.

Let $\mathscr H_*$ be the product for $k\neq i$ of the complete flag
varieties of the vector spaces $V_k$. Let $\mathscr H'_i$ be the
set of all filtrations
$$0=V'_{i,0}\subseteq V'_{i,1}\subseteq\cdots\subseteq
V'_{i,r-1}\subseteq V'_{i,r}=V'_i$$
such that $\dim V'_{i,s}/V'_{i,s-1}=m_s$ for each $s\in[1,r]$.
Certainly, $\mathscr H'$ is isomorphic to the product
$\mathscr H'_i\times\mathscr H_*$. Our integral can then be computed
with the help of the Fubini theorem:
$$\langle\delta_M,u\rangle=\int_{\mathscr H'}\phi'=
\int_{\mathscr H_*}\Biggl(\int_{\mathscr H'_i}\phi'\Biggr).$$

The right hand side of (\ref{eq:GenForm}) can be computed by a similar
convolution product. Let $\mathscr H''_i$ be the set of all filtrations
$$0=V''_{i,0}\subseteq V''_{i,1}\subseteq\cdots\subseteq
V''_{i,r-1}\subseteq V''_{i,r}=V''_i$$
such that $\dim V''_{i,s}/V''_{i,s-1}=-a_{ij_s}-m_s$ for each
$s\in[1,r]$. Then
$$\langle\delta_M,T_i(u)\rangle=
\int_{\mathscr H_*}\Biggl(\int_{\mathscr H''_i}\phi''\Biggr),$$
where $\phi''$ is a function given by the product of the
$\kappa\bigl((-1)^{m_s}f_{i,j_s,-a_{ij_s}-m_s}\bigr)$.

Equation (\ref{eq:GenForm}) can thus be written
$$\int_{\mathscr H_*}\Biggl(\int_{\mathscr H'_i}\phi'\Biggr)=
\int_{\mathscr H_*}\Biggl(\int_{\mathscr H''_i}\phi''\Biggr).$$
Therefore to prove (\ref{eq:GenForm}), one only has to show
\begin{equation}
\label{eq:PartGenForm}
\int_{\mathscr H'_i}\phi'=\int_{\mathscr H''_i}\phi'',
\end{equation}
where each side depends on the choice of a point in $\mathscr H_*$.

\subsection{}
\label{ss:SplittingArg}
We keep the notations of section~\ref{ss:PartInt}.

For each $k\neq i$, let $A_k$ be the vector space with basis
$\{h\in H\mid(s(h),t(h))=(i,k)\}$. Let $\widetilde V=\widetilde M_i$;
in~other~words,
$$\widetilde V=\bigoplus_{\substack{h\in H\\s(h)=i}}
V_{t(h)}=\bigoplus_{k\neq i}V_k\otimes_{\mathbb C}A_k.$$
We also set $x=M_{\out(i)}M_{\iin(i)}$.

We have chosen a point in $\mathscr H_*$, that is, a filtration in
$V_k$ for each $k\neq i$:
$$0=V_{k,0}\subseteq V_{k,1}\subseteq\cdots\subseteq V_{k,r-1}\subseteq
V_{k,r}=V_k,$$
such that $\dim V_{k,s}/V_{k,s-1}=1$ if $k=j_s$ and $V_{k,s}=V_{k,s-1}$
if $k\neq j_s$. This induces a filtration
\begin{equation}
\label{eq:FiltTildeV}
0=\widetilde V_0\subseteq\widetilde V_1\subseteq\cdots\subseteq
\widetilde V_{r-1}\subseteq\widetilde V_r=\widetilde V,
\end{equation}
namely
$$\widetilde V_s=\bigoplus_{k\neq i}V_{k,s}\otimes_{\mathbb C}A_k.$$

We set $n=\dim\widetilde V$. For each $s\in[1,r]$, we set
$p_s=1+\dim\widetilde V_{s-1}$ and $q_s=\dim\widetilde V_s$. We thus
have a partition $[1,n]=[p_1,q_1]\sqcup\cdots\sqcup[p_r,q_r]$ as a
union of disjoint intervals.

Let us first assume that $x$ does not leave stable the filtration
(\ref{eq:FiltTildeV}). In this case, both sides of (\ref{eq:PartGenForm})
vanish, so (\ref{eq:PartGenForm}) holds true.

Let us now assume that $x$ leaves stable each $\widetilde V_s$.
Noting that $x$ induces a nilpotent endomorphism on each quotient
$\widetilde V_s/\widetilde V_{s-1}$, we can find a basis
$\mathbf e=(e_1,\ldots,e_n)$ of $\widetilde V$ in which the matrix of
$x$ is strictly upper triangular and such that $(e_1,\ldots,e_{q_s})$
is a basis of $\widetilde V_s$, for each $s\in[1,r]$.

We then consider the graph (\ref{eq:StarDiag}). The objects attached
to this graph will be underlined: the enveloping algebra is
$\underline{\mathbf U}^+$ and is generated by elements $\underline e{}_i$
with $i\in\{0,1,\ldots,n\}$; the completed preprojective algebra
is $\underline\Pi$.

We set $\underline M{}_0=V'_i$ and $\underline M{}_j=\mathbb Ce_j$ for
all $j\in[1,n]$; thus $M_i=\underline M{}_0$ and
$\widetilde M_i=\underline M{}_1\oplus\cdot\oplus\underline M{}_n$.
Writing the linear maps $M_{\iin(i)}$ and $M_{\out(i)}$ as block
matrices with respect to this decomposition of $\widetilde M_i$, we
get maps $\underline M{}_j\to\underline M{}_0$ and
$\underline M{}_0\to\underline M{}_j$. These maps satisfy the
preprojective relations of $\underline\Pi$, thanks to
Equation~(\ref{eq:MinMout}) and to the fact that the matrix of $x$
in $\mathbf e$ is strictly upper triangular. We therefore get a
$\underline\Pi$-module $\underline M$ such that
$\underline M{}_{\out(0)}=M_{\out(i)}$ and
$\underline M{}_{\iin(0)}=M_{\iin(i)}$.

With these notations, the integral on the left hand side of
(\ref{eq:PartGenForm}) computes
$$\langle\delta_{\underline M},\underline g{}'_r\cdots\underline
g{}'_1\rangle,\qquad\text{where}\quad\underline
g{}'_s=\sum_{p+q=m_s}(-1)^p\;\frac{\underline
e{}_0^p}{p!}\;\underline e{}_{q_s}\cdots\underline
e{}_{p_s}\;\frac{\underline e{}_0^q}{q!}.$$

Replacing $M$ by $\Sigma_i^*M$ amounts to replace $M_i$ by
$\coker M_{\out(i)}$ without changing $\widetilde M_i$ nor $x$.
The analogous of $\underline M$ for $\Sigma_i^*M$ is therefore
$\Sigma_0^*\underline M$. Thus the right hand side of
(\ref{eq:PartGenForm}) is equal to
$$\langle\delta_{\Sigma_0^*\underline M},\underline
g{}''_r\cdots\underline g{}''_1\rangle,\qquad\text{where}\quad\underline
g{}''_s=(-1)^{m_s}\sum_{p+q=-a_{ij_s}-m_s}(-1)^p\;\frac{\underline
e{}_0^p}{p!}\;\underline e{}_{q_s}\cdots\underline
e{}_{p_s}\;\frac{\underline e{}_0^q}{q!}.$$

Let $D=\ad(\underline e{}_0)$, a derivation of the algebra
$\underline{\mathbf U}^+$. The Serre relations say that
$D^2(\underline e{}_j)=0$ for any $j\in[1,n]$. Using Leibniz's
formula for the $m_s$-th derivative of a product, we obtain
$$\underline g{}'_s=\frac{(-1)^{m_s}}{m_s!}D^{m_s}(\underline
e{}_{q_s}\cdots\underline e{}_{p_s})=\sum_{\substack{\strut
J\subseteq[p_s,q_s]\\|J|=m_s}}\underline
v{}_{q_s}^J\cdots\underline v{}_{p_s}^J,$$
where
$$\underline v{}_j^J=\begin{cases}
[\underline e{}_j,\underline e{}_0]=\underline f{}_{0,j,1}&
\text{if $j\in J$,}\\
\underline e{}_j=\underline f{}_{0,j,0}&\text{if $j\notin J$.}
\end{cases}$$
With the same notation,
$$\underline g{}''_s=\frac{(-1)^{a_{ij_s}}}{(-a_{ij_s}-m_s)!}
D^{-a_{ij_s}-m_s}(\underline e{}_{q_s}\cdots\underline e{}_{p_s})
=(-1)^{m_s}\sum_{\substack{\strut
K\subseteq[p_s,q_s]\\|K|=-a_{ij_s}-m_s}}\underline
v{}_{q_s}^K\cdots\underline v{}_{p_s}^K.$$

Matching the summand indexed by $J$ in the expansion of $\underline
g{}'_s$ with the summand indexed by $K=[p_s,q_s]\setminus J$ in the
expansion of $\underline g{}''_s$, we get $\underline
g{}''_s=T_0(\underline g{}'_s)$, so
$$\langle\delta_{\Sigma_0^*\underline M},\underline
g{}''_r\cdots\underline g{}''_1\rangle=
\langle\delta_{\Sigma_0^*\underline M},T_0(\underline
g{}'_r\cdots\underline g{}'_1)\rangle.$$
We then see that (\ref{eq:PartGenForm}) is simply the result of
section~\ref{ss:StarDiag} applied to $\underline M$ and to
$\underline u=\underline g{}'_r\cdots\underline g{}'_1$.

This last argument finishes the proof of (\ref{eq:PartGenForm}),
hence of (\ref{eq:GenForm}).

Pierre Baumann\\
Institut de Recherche Math\'ematique Avanc\'ee\\
Universit\'e de Strasbourg et CNRS\\
7 rue René Descartes\\
67084 Strasbourg Cedex\\
France\\[5pt]
\texttt{p.baumann@unistra.fr}
\end{document}